\documentclass[12pt]{amsart}
\usepackage{amsmath}

\def\ol{\overline}
\def\e{\epsilon}
\def\lf{\left}
\def\ri{\right}

\def\wt{\widetilde}

\def\p{\partial}

\newcommand\C{{\mathbb C}}

\def\jbar{{\bar\jmath}}

\def\K{K\"ahler }
\def\KR{K\"ahler-Ricci }
\def\KRF{K\"ahler-Ricci flow }
\def\KRS{K\"ahler-Ricci soliton }
\def\KRS{K\"ahler-Ricci soliton }

\def\be{\begin{equation}}
\def\ee{\end{equation}}
\def\ol{\overline}
\def\lf{\left}
\def\ri{\right}
\def\aut{\text{\rm Aut}}

\def\e{\epsilon}

\def\ijb{{i\jbar}}

\def\wt{\widetilde}

\def\p{\partial}

\def\C{\Bbb C}

\def\cn{\Bbb C^n}
\def\wh{\widehat}

\def\wt{\widetilde}

\def\p{\partial}

\def\C{\Bbb C}
 
\def\KRF{K\"ahler-Ricci flow }

\newtheorem{thm}{Theorem}[section]

\newtheorem{lem}{Lemma}[section]

\newtheorem{cor}{Corollary}[section]
\theoremstyle{definition}
\newtheorem{defn}{Definition}[section]
\theoremstyle{remark}
\newtheorem{rem}{Remark}[section]
\numberwithin{equation}{section}
 \begin{document}
\author{Albert Chau$^1$}
\thanks{$^1$Research
partially supported by NSERC grant no. \# 327637-06}

\address{Waterloo University, Department of Pure Mathematics,
  200 University avenue, Waterloo, ON N2L 3G1, CANADA}
\email{a3chau@math.uwaterloo.ca}

\author{Luen-Fai Tam$^2$}
\date{October,  2006; revised August, 2007}
\thanks{$^2$Research
partially supported by Earmarked Grant of Hong Kong \#CUHK403005}
\address{Department of Mathematics, The Chinese University of Hong Kong,
Shatin, Hong Kong, China.} \email{lftam@math.cuhk.edu.hk}

\title{On the Steinness of a class of K\"ahler manifolds}

\maketitle \markboth{Albert Chau and Luen-Fai Tam} {Steinness of a class of K\"ahler manifolds}

\section{Introduction}
Let $(M^n, g_0)$ be a complete non-compact \K manifold with
complex dimension $n$ and with bounded nonnegative holomorphic
bisectional curvature. Let $R$ be the scalar curvature and define
$$k(x, r):=\frac{1}{V_x(r)}\int_{B_x(r)}R dV.$$ In \cite{CT3}, it was
proved by the authors that if $M$ has maximum volume growth, then $M$ is
biholomorphic to
$\cn$.  There, the authors used a result of Ni in \cite{Ni2} (see also \cite{CTZ,CZ3}) which
states that the condition of maximum volume growth on $M$ implies that
\begin{equation}\label{quadraticdecay}
k(x,r)\le \frac{C}{1+r^2}
\end{equation}
for some $C$ for all $x$ and $r$. In \cite{CT4}, the authors
proved that condition (\ref{quadraticdecay}) implies that $M$ is
holomorphically covered by $\cn$, without assuming the maximum
volume growth condition. The proof is obtained by  studying the
\KRF:
\begin{equation}\label{s1e0}
   \frac{d g_\ijb}{dt}=-R_\ijb
\end{equation}
with initial data $g_0$. It is well-known by \cite{Sh2} that if
the scalar curvature decays linearly in the average sense:
\begin{equation}\label{lineardecay1}
k(x,r)\le C/(1+r)
\end{equation}
for some constant $C$ for all $x$ and $r$ then (\ref{s1e0})
has long time solution with uniformly bounded curvature. By
the results in \cite{CZ3,NT2}, the linear decay condition
(\ref{lineardecay1}) is true in most case, at least for
constant $C$ which may depend on $x$.

 In this paper, we will prove the
following:

\begin{thm}\label{s1t1}
Let $(M^n, g_0)$ be a complete non-compact \K manifold with
bounded non-negative holomorphic bisectional curvature.
Suppose the scalar curvature of $g_0$ satisfies the linear
decay condition (\ref{lineardecay1}). Then $M$ is
holomorphically covered by a pseudoconvex domain in $\C^n$
which is homeomorphic to $\mathbb{R}^{2n}$.  Moreover, if
$M$ has positive bisectional curvature and is simply
connected at infinity, then M is biholomorphic to a
pseudoconvex domain in $\cn$ which is homeomorphic to
$\mathbb{R}^{2n}$, and in particular, $M$ is Stein.

\begin{rem}\label{s0r1}
By a result of Yau \cite{Y1}, the pseudoconvex domain in Theorem
\ref{s1t1}  has infinite Euclidean volume. The authors would like
to thank Shing Tung Yau for providing this information.
\end{rem}
\end{thm}

If we assume that $k(r)=\frac{C}{1+r^{1+\epsilon}}$, for $\epsilon
>0$, the result that $M$ is biholomorphic to a pseudoconvex domain
was proved by Shi \cite{Sh2} under the additional assumption that
$(M, g)$ has positive sectional curvature. Note that if $M$ has
positive sectional curvature, then it is well-known that $M$ is
diffeomorphic to $\mathbb{R}^{2n}$ by \cite{GM}, and is  Stein by
  \cite{GW76}.  Under the same decay condition and assuming maximum
volume growth, similar results were obtained by Chen-Zhu
\cite{Chen}. All these works are before   \cite{CT3,CT4}.

As in the above mentioned works, our proof of Theorem
\ref{s1t1} is based on the  \KR flow (\ref{s1e0}). In fact,
Theorem \ref{s1t1} will be proved as a consequence of the
following more general:

\begin{thm}\label{main} Let $M^n$ be a complex noncompact manifold.
 Suppose there exist a
 sequence of complete K\"ahler metrics $g_i$, for $i\ge1$, on $M$ such
 that
 \begin{enumerate}
     \item  [\textbf{ (a1)} ] $cg_{i}\le g_{i+1}\le g_i$  for some
     $1>c>0$ for all $i$.

     \item  [\textbf{ (a2)} ] $|Rm(g_i)|+|\nabla Rm(g_i)|\le c'$ for some $c'$
         on $ B_i(p,r_0)$ for some $p\in M$ and $r_0>0$ for all $i$
         where  $B_i(p,r_0)$ is the geodesic ball with respect to
$g_i$.

 \item  [\textbf{ (a3)} ]  $g_i$ is contracting in the
following sense:   For any $\e$,  for any $i$, there exists $i'>i$
with
$$
g_{i'}\le \e g_i
$$
in $B_i(p, r_0)$.
\end{enumerate}
  Then  $M$ is covered by a pseudoconvex
domain in $\cn$ which is homeomorphic to $\mathbb{R}^{2n}$.
\end{thm}
   To prove Theorem \ref{s1t1},   the solution
$g(t)$ to the \KR flow on $M$ will be used to produce a
sequence of \K metrics $g_i$ satisfying the hypothesis of
the Theorem \ref{main}.  The main steps in proving Theorem
\ref{main} can be sketched as follows. The idea is to
consider the sequence of holomorphic normal ``coordinate
charts'' around some $p \in M$ corresponding to the
sequence $g_i$. We then use this sequence of charts,
together with a gluing technique as in \cite{Sh2}, to build
a map from an open set in $\C^n$ onto $M$.  In general
however, these charts will only be locally biholomorphic,
and to build such a map one generally needs to control the
sets around $p$ on which these charts are injective
\footnote{In \cite{Sh2}, positive sectional curvature was
used to produce a sequence of strictly convex domains
around $p$ exhausting $M$, which were then used to control
the injectivity of the charts. In \cite{Chen}, maximal
volume growth was used to control the injectivity radius
under the \KR flow}.  We will not assume any control these
sets. Instead, we will control the sets where the
corresponding coordinate transition functions are injective
using a method developed by the authors in \cite{CT4}. Once
the transition functions are established, we then following
similar techniques in \cite{Sh2} and \cite{CT4}, to build a
covering map from an open set $\Omega$ in $\C^n$ onto $M$.
By its construction, $\Omega$ will be shown to be
pseudoconvex and homeomorphic to $\mathbb{R}^{2n}$.

\section{holomorphic coordinate ``covering" charts}

Let $M, g_i, p, r_0$ be as in Theorem \ref{main}.

\begin{lem} \label{s1-L1}
  There exists $r>0$, and  a family of holomorphic maps
\begin{equation}\label{s1e2}
\Phi_i:D(r)\to M
\end{equation}
 for all $i\geq 1$ with the following
  properties:
  \begin{enumerate}
\item [(i)] $\Phi_i$ is a local biholomorphism from $D(r) \subset
\C^n$ onto its image,

\item[(ii)] $\Phi_i(0)=p$,

\item[(iii)] $\Phi_i^*(g_i))(0)=g_e$,

\item[(iv)] $\frac{1}{C}g_{e}\leq\Phi_i^*(g_i) \leq  Cg_{e}$ in
$D(r)$,
\end{enumerate}
 where $g_e$ is the standard metric on $\cn$, $C$ is a constant independent of
$t$ and $p$.
\end{lem}
\begin{proof} Using condition \textbf{(a2)} and by considering the pullback metric under the
exponential map within the conjugate locus, one can   apply
Proposition 2.1 in \cite{CT3} to obtain the results.
\end{proof}


\begin{cor}\label{cor-to-s1-L1}
 $B_i(C^{-1}\rho)\subset\Phi_i(D(\rho))\subset
B_i(C\rho)$ for some $C>0$ for all   $ 0<\rho<r$ and $i\ge1$.
\end{cor}

The following two lemmas  are   from \cite[Lemmas 3.2 and
3.3]{CT4}.

\begin{lem}\label{s1-L2}
For any $0<\rho \le r$,  where $r$ is as in Lemma \ref{s1-L1},
there exists $0<\rho_1<r_0$ , independent of $i$, satisfying the
following:

 \begin{itemize}
\item [(i)]  For any $q\in B_{i}(p,\rho_1)$, there is $z\in
D(\frac\rho8)$ such that $\Phi_i(z)=q$.

\item[(ii)] For any  $q\in B_{i}(p,\rho_1)$,  $z\in D(\frac\rho8)$
with $\Phi_i(z)=q$, and any smooth curve $\gamma$ in $M$ with
$\gamma(0)=q$ such that $L_i(\gamma)<\rho_1$, there is a unique
lift $\wt\gamma$ of $\gamma$ by $\Phi_i$ so that $\wt\gamma(0)=z$
and $\wt\gamma\subset D(\frac\rho2 )$.
\end{itemize}
\end{lem}

\begin{lem}\label{s1-L3} Fix $i\ge1$.  Let $0<\rho\le r$
be given and let $\rho_1$ be as in Lemma \ref{s1-L2}.  Given any
$\e>0$,
 there exists $\delta>0$, which may depend on $i$,
 satisfying the following properties:

 Let  $\gamma(\tau)$, $\beta(\tau)$, $\tau\in [0,1]$ be  smooth
curves from $q\in B_{i}(p,\rho_1  )$ with length less than $\rho_1
 $ with respect to $g_i$ and let $z_0\in D(\frac18\rho )$
with $\Phi_i(z_0)=q$. Let $\wt\gamma$, $\wt\beta$ be the
liftings from $z_0$ of $\gamma$ and $\beta$ as described in
Lemma \ref{s1-L2}.
 Suppose  $d_i(\gamma(\tau),\beta(\tau))<\delta$ for all $\tau\in [0,1]$,
 then $d_e(\wt\gamma(1),\wt\beta(1))<\e$.
 Here   $d_i$ is the distance in $g_i$ and $d_e$ is the Euclidean distance.
\end{lem}
\begin{cor}\label{cor-to-L3} Let $0<\rho\le r$
be given and let $\rho_1$ be as in Lemma \ref{s1-L3}. Let
$\gamma:[0,1]\times[0,1]\to M$ be smooth homotopy such that
 \begin{enumerate}
    \item[(a)] $\gamma(0,\tau)=q_1$ and $\gamma(1,\tau)=q_2$ for all
    $\tau$.
    \item[(b)] $q_1\in B_i(p,\rho_1 )$ and $\Phi_i(z_0)=q_1$
    for some $z_0\in D(\frac18\rho r(t))$.
    \item[(c)] For all $0\le \tau\le 1$, the length of
    $\gamma(\cdot,\tau)$ is less than $\rho_1$.
\end{enumerate}
For all $\tau$, let $\wt  \gamma_\tau$ be the lift of
$\gamma(\cdot,\tau)$ as in Lemma \ref{s1-L2} from   $z_0$.
Then $\wt\gamma_\tau(1)=\wt\gamma_0(1)$ for all $\tau$.
\end{cor}
\begin{proof}  By Lemma \ref{s1-L2}, $\wt\gamma_\tau(1)\in D(\frac12\rho )$ for
all $\tau$.   Let $\e>0$ be such that $\Phi_{i}$ is
injective on $D(w,\e)$ for all $w\in D(\frac12\rho  )$. Let
$\delta>0$ be as in Lemma \ref{s1-L3}.  Let $m$ be large
enough, such that
$d_{i}(\gamma(s,j/m),\gamma(s,(j+1)/m))<\delta$ for all $s$
for $0\le j\le m-1$. By Lemma \ref{s1-L3}, we have
$|\wt\gamma_{j/m}(1)-\wt\gamma_{(j+1)/m}(1)|<\e$ for $0\le
j\le m-1$. Since $\Phi_{i}\circ \wt\gamma_\tau(1)=q_2$,
 and $\Phi_{i}$ is injective in $D(\wt\gamma_\tau(1),\e)$, we have
 $$
 \wt\gamma_{j/m}(1)=\wt\gamma_{(j+1)/m}(1).
 $$
 From this the Corollary follows.
\end{proof}

\section{holomorphic transition functions}

Let $M, g_i, p, r_0$ be as in Theorem \ref{main}.
\begin{lem}\label{s2-L1} Let $r$ be as in Lemma \ref{s1-L2}.
There exists $r>\rho>0$ such that for every $i\geq 1$, where there
is a map $F_{i+1}$ from $ D(\rho  )$ to $ D(r )$ such that
$\Phi_i=\Phi_{i+1}\circ F_{i+1}$ on $ D(\rho  )$. Moreover,
$\Phi_i(D(\rho))\subset B_i(p,r_0)$ where $r_0$ is the constant in
\textbf{(a3)}.
\end{lem}
\begin{proof}
In Lemma \ref{s1-L2}, let $\rho=r$ and let $\rho_1$ be as in the
conclusion of the Lemma.  Note that $\rho_1$ is independent of
$i$.  Now let $0<K<1$ be a constant to be determined.  For any
$z\in D(K \rho_1 )$, let $\gamma^*(\tau)$, $0\le \tau\le 1$, be
the line segment from 0 to $z$, and let
$\gamma=\Phi_i\circ\gamma^* $.  By \textbf{(a1)}and Lemma
\ref{s1-L1}, we see that there is a constant $C_1>0$ independent
of $i$   such that
\begin{equation}
  L_{i+1}(\gamma)\le L_i(\gamma)< C_1K \rho_1.
    \end{equation}
Now choose $K$  so that $C_1 K<1$ and let $\rho = K \rho_1$. Since
$\gamma(0)=p$,    by Lemma \ref{s1-L2}, there is a unique lift
  $\wt\gamma$
  of $\gamma$ by $\Phi_{i+1}$ so that $\wt\gamma(0)=0$ and
  $\wt\gamma\subset D(\frac12 r)$. We  define
  $F_{i+1}(z)=\wt\gamma(1)$. $F_{i+1}$ is then a well-defined map from
  $  D(\rho   )$ to  $D(r )$
   by  the uniqueness of  the lifting.  Also, by construction we have
   $\Phi_i=\Phi_{i+1}\circ
F_{i+1}$ on $ D(\rho)$. By choosing a smaller $\rho$, we
also have $\Phi_i(D(\rho))\subset B_i(p,r_0)$.
  This completes the proof of the Lemma.
  \end{proof}

  \begin{lem}\label{transitionmaps} Let $\rho$ be as in Lemma \ref{s2-L1}.
  Then for any $i\geq 1$, the map $F_{i+1}$ satisfies the
  following:

  \begin{enumerate}
    \item [(a)] $F_{i+1} (0)=0$
    \item [(b)] $F_{i+1}$ is a local  biholomorphism.
    \item [(c)]
    $$
    b_1|v|\le |F'_{i+1}(0)v|\le b_2|v|
    $$
    for some $0<b_1\le b_2\leq1$ independent of $i$,
    and for all vectors $v\in
    \cn$, where $F'$ is the Jacobian of $F$.
       \item[(d)] There exist  $\rho_1$ and $\rho_2$
       independent of $i$, each in $(0, \rho)$, such that
       $$F_{i+1}(D(\rho_1 )) \subset D(\rho_2 ),$$
       and $F^{-1}_{i+1}$ exists on $D(\rho_2 )$.
 \end{enumerate}
  \end{lem}
\begin{proof}
(a) follows from the definition of $F_{i+1}$.

(b) can be proved as in the proof of Lemma 3.4 part (b) in
\cite{CT4}, using Lemmas \ref{s1-L3} and \ref{s2-L1}.

(c) follows from the \textbf{(a1)} and Lemma \ref{s1-L1}.

 (d)  follows from the   proof of part (d) of Lemma 3.4 in
 \cite{CT4}.
\end{proof}
\begin{cor}\label{s2-C1}
Let $\rho_1$ be as in Lemma \ref{transitionmaps}. Then for any
$i\geq 1$, $F_{i+1}(D(\rho_1))$ is Runge in $\C^n$.
\end{cor}
\begin{proof}
Let $i\geq 1$ be given.   Then   given any holomorphic function
$f$ on $F_{i+1}(D(\rho_1)) \subset \C^n$, we must show that $f$
can be approximated by entire functions on $\C^n$ uniformly on
compact subsets of $F_{i+1}(D(\rho_1)) $.   Consider the
holomorphic function $ f\circ F_{i+1}$ defined on $D(\rho_1 )$.
Since $D(\rho_1  )$ is just a ball in $\C^n$, $ f\circ F_{i+1}$
can be approximated uniformly on compact subsets of $D(\rho_1 )$
by entire functions $h$ on $\C^n$. By part (d) of Lemma
\ref{transitionmaps},   $ h\circ F_{i+1}^{-1} $ are defined on
$D(\rho_2)$ and holomorphic. We see that these approximate $f$
uniformly on compact subsets of $F_{i+1}(D(\rho_1 ))\subset
D(\rho_2 )$.  Finally, as $D(\rho_2 )$ is just a ball in $\C^n$,
we see that the functions $h\circ F_{i+1}^{-1}$ can themselves be
approximated uniformly on compact subsets of $D(\rho_2 )$ by
entire functions. Thus by part (d) of Lemma \ref{transitionmaps}
$f$ can be uniformly approximated on compact subsets of
$F_{i+1}(D(\rho_1))$ by entire functions. This completes the proof
of the Corollary.
\end{proof}
\begin{cor}\label{s2-C2}
Let $\rho_1$ be as in Lemma \ref{transitionmaps}. Then for any
$i\geq 1$, $F_{i+1}$ can be approximated uniformly on compact
subsets of $D(\rho_1  )$ by elements of $ {\rm Aut}(\C^n)$.
\end{cor}

\begin{proof}
This follows from Corollary \ref{s2-C1} and Theorem 2.1 in
\cite{A}.
\end{proof}

 \section{construction of a map onto $M$}

Let $M, g_i, p, r_0$ be as in Theorem \ref{main}. We begin with
the following lemma on the transition functions $F_i$ which
basically says that the maps are contracting.

\begin{lem}\label{fixl1}
Let $\rho_1$ be as in Lemma \ref{transitionmaps}.  Then there
exists positive constants $\rho_3$ and $C$  with $C>1$,
  $C\rho_3 < \rho_1$,   such that for every $i$ and $k\ge 1$
\begin{equation}\label{fixe0}
F_{i+k}\circ \cdot \cdot \circ F_{i+1}(D(\rho_3 )) \subset
D(C\rho_3 ).
\end{equation}
\end{lem}
\begin{proof}
Let $C_1>1$ be the constant in property (iv) of $\Phi_i$ in Lemma
\ref{s1-L1}  and $\rho_1$ as in Lemma \ref{transitionmaps}. Let
$C=C_1^2 $ and $\rho_3=\rho_1/(2C)$. For any $i$, we want to prove
that
$$
F_{i+1}(D(\rho_3 ))\subset  D(C\rho_3).
$$
Let $$A=\{\eta\in(0, \rho_3 ]|\ F_{i+1}(D(\eta))\subset D(C\rho_3
) \}
$$
Since $F_{i+1}(0)=0$ and $F_{i+1}$ is a local biholomorphism, it
is easy to see that $A$ is nonempty   in $(0, \rho_3  ]$. Let
$r\in A$ and $z\in D(\eta)$. Since $C\rho_3<\rho_1$, $\Phi_i(tz)$
and $\Phi_{i+1}\circ F_{i+1}(tz)$, $0\le t\le 1$, are defined and
are equal by Lemmas \ref{s2-L1}, \ref{transitionmaps}. By Lemma
\ref{s1-L1}, we have
\begin{equation}\label{s4f1e}
\begin{split}
||\frac{d}{dt}F_{i+1}(tz)||_{g_e}&\le
C_1||\frac{d}{dt}F_{i+1}(tz)||_{\Phi_{i+1}^*(g(i+1))}\\
&=C_1||\frac{d}{dt}\Phi_{i+1}\circ F_{i+1}(tz)||_{ g(i+1) }\\
&\le C_1||\frac{d}{dt}\Phi_{i}(tz)||_{ g(i) }\\
&= C_1||\frac{d}{dt}(tz)||_{\Phi_i^*(g(i))}\\
&\le C_1^2|z|\\
&< C_1^2 \eta\\
&\le C\rho_3,
\end{split}
\end{equation}
where we have used   \textbf{(a1)}. From this, it is easy to see
that $A$ is both open and closed in $(0, \rho_3  ]$ and
$F_{i+1}(D(\rho_3 ))\subset D(C\rho_3 )$.

Suppose $k> 1$ such that
\begin{equation}
 F_{i+l}\circ \cdot \cdot \circ F_{i+1}(D(\rho_3  ))\subset D(C\rho_3  )
\end{equation}
for all $1\le l< k$. As before, let
$$B=\{\eta\in(0, \rho_3  ]|\
F_{i+k}\circ \cdot \cdot \circ F_{i+1}(D(\eta))\subset D(C\rho_3
 )   \}.
$$
Again, $B$ is nonempty   in $(0, \rho_3  ]$. Suppose $r\in B$ and
$z\in D(r)$. Then $\Phi_i(tz)$ and $\Phi_{i+k}\circ F_{i+k}\circ
\cdot \cdot \circ F_{i+1}(tz)$,
 $0\le t\le 1$,
are well-defined and equal. As before, we can prove that $B$ is
open and closed in $(0, \rho_3 ]$ and
\begin{equation}
 F_{i+k}\circ \cdot \cdot \circ F_{i+1}(D(\rho_3  ))\subset
 D(C\rho_3  ).
 \end{equation}
  This completes the
proof of the Lemma.
\end{proof}
\begin{rem}\label{s5r1} For later use, we will assume that
$\rho_3<\frac18 r$, where $r$ is as in Lemma \ref{s1-L1}.
\end{rem}

\begin{lem}\label{fixl22} Let $\rho_3$ as in Lemma \ref{fixl1}
There exists a positive increasing sequence $n_i$ for $i\geq 1$
such that $n_1=1$ and
\begin{equation}\label{contracting1}
F_{n_{i+1}}\circ \cdot \cdot F_{n_i +2}\circ F_{n_i
+1}(D(\rho_3))\subset D(\frac {\rho_3}2)
\end{equation}
for every $i$.
\end{lem}
\begin{proof} Let $C$ be the constant in Lemma \ref{fixl1}.
Let $n_1=1$. By Lemma \ref{fixl1}, for all $k$, $
F_{n_1+k}\circ\cdots\circ F_{n_1+1}$ is defined in
$D(\rho_3)$ for all $k\ge1$. As in the
proof of Lemma \ref{fixl1}, for all $z\in D(\rho_3)$, and
$0\le t\le 1$:
\begin{equation}
\begin{split}
||\frac{d}{dt}F_{n_1+k}\circ\cdots\circ
F_{n_1+1}(tz)||_{g_e}&\le
C_1||\frac{d}{dt}F_{n_1+k}\circ\cdots\circ
F_{n_1+1}(tz)||_{\Phi_{k+n_1}^*(g_{n_1+1}) }\\
&=C_1||\frac{d}{dt}\Phi_{k+n_1}\circ F_{n_1+k}\circ\cdots
\circ F_{n_1+1}(tz)||_{ g_{k+n_1} }\\
&= C_1||\frac{d}{dt}\Phi_{n_1}(tz)||_{ g_{k+n_1}}
\end{split}
\end{equation}
where $C_1$ is as in the proof Lemma \ref{fixl1}. Since
$\Phi_{n_1}(D(\rho_3))\subset B_{n_1}(p,r_0)$ by the choice
of $\rho$ in Lemma \ref{s2-L1}, by \textbf{(a3)} and (iv)
in Lemma \ref{s1-L1},  we can find $n_2>n_1$ such that
$$
F_{n_2}\circ  \cdots\circ F_{n_1+1}(D(\rho_3))\subset
D(\frac12\rho_3).
$$
Similarly, one can choose $n_3,n_4,\dots$ inductively which
satisfy the conclusion of the lemma.

\end{proof}

We now want to construct an appropriate sequence
$\tilde{F}_j \in {\rm Aut}(\C^n)$ which will approximate
the sequence $F_j$ for $j\geq 2$.

Let $n_i$ be as in Lemma \ref{fixl22}. By Lemmas
\ref{fixl1}, \ref{fixl22} and Corollary \ref{s2-C2}, we can
find  $\tilde F_2,\dots,\tilde F_{n_2}$ in $\aut(\cn)$ such
that
\begin{equation}\label{approx1}
    \tilde F_{k+1}\circ\cdots \circ\tilde F_2(D(\rho_3))\subset
    D(\rho_1)
\end{equation}
for $2\le k\le n_2$ and
\begin{equation}\label{approx2}
    \tilde F_{n_2}\circ\cdots \circ\tilde F_2(D(\rho_3))\subset
    D( \rho_3).
\end{equation}

Since $\Phi_{n_i}$ is a local biholomorphism, we have
\begin{equation}\label{approx3}
    ||D\Phi_{n_2}\circ\tilde F_{n_2}\circ\cdots
    \circ\tilde F_2(z)(v)||_{g_1}\ge b_2>0
\end{equation}
for all  $z\in D(\rho_3)$ and unit vectors $v\in \cn$.


Let $S_2=(\tilde F_{n_2}\circ \cdots \tilde
F_2)^{-1}(D(\rho_3))$. Use Lemmas \ref{fixl1}, \ref{fixl22}
and Corollary \ref{s2-C2} again, we can find $n_3>n_2$ and
$\tilde F_{n_2+1},\dots,\tilde F_{n_3}$ in $\aut(\cn)$ such
that
\begin{equation}\label{approx5}
    \tilde F_{n_3}\circ\cdots \circ\tilde F_{n_2+1}(D(\rho_3))\subset
    D(\rho_3).
\end{equation}
Since
$$ \Phi_{n_2}  =\Phi_{n_3}\circ F_{n_3}\circ\cdots\circ
F_{n_2+1}
$$
on $D(\rho_3)$, we may choose $\tilde F_{n_2+1},\dots,\tilde
F_{n_3}$ such that they also satisfy:
\begin{equation}\label{approx5}
   d_{g_1}(\Phi_{n_{3}} \circ \tilde{F}_{n_{3}} \circ
\cdots\circ\tilde F_{n_2+1}\circ \tilde F_{n_2}\circ\cdots
\circ\tilde{F}_{2} (z), \Phi_{n_2} \circ \tilde{F}_{n_2}\circ
\cdots v \circ \tilde{F}_{2}(z))\leq \frac{1}{2^{2}}
\end{equation}
for $z\in S_2$ and
\begin{equation}\label{approx6}
\begin{split}
||D(\Phi_{n_3}& \circ \tilde{F}_{n_3}\cdots
\circ\tilde{F}_{n_2+1}\circ\tilde{F}_{n_2}\circ \cdots
\circ
\tilde{F}_{2})(z)(v)||_{g_1}\\
&-||D(\Phi_{n_2} \circ \tilde{F}_{n_2}\circ \cdots  \circ
\tilde{F}_{2})(z)(v)||_{g_1} \leq \frac{ b_2 }{2^{2}}
\end{split}
\end{equation}
for all $z\in S_2$, and for all unit vector in $\cn$. Let
$0<b_3<b_2$ be such that
\begin{equation}\label{approx3}
    ||D\Phi_{n_3}\circ\tilde F_{n_3}\circ\cdots
    \circ\tilde F_2(z)(v)||_{g_1}\ge b_3>0
\end{equation}
for all  $z\in S_2$ and unit vector $v\in \cn$. Let
$S_3=(\tilde F_{n_3}\circ \cdots \tilde
F_2)^{-1}(D(\rho_3)$. Inductively in this way, we have the
following:

\begin{lem}\label{fixc2}
There exist $\tilde F_2,\dots, \tilde F_j,\dots$ in
$\aut(\cn)$, such that  the following conditions are
satisfied for all $i\geq2 $:

\begin{equation}\label{fixe8}
\tilde{F}_{n_{i+1}}\circ \cdots  \circ \tilde{F}_{n_i
+1}(D(\rho_3 ))\subset D(\rho_3  )).
\end{equation}

\begin{equation}\label{approxfix}
d_{g_1}(\Phi_{n_{i+1}} \circ \tilde{F}_{n_{i+1}} \cdot
\cdot \circ\tilde{F}_{n_{i}+1}\circ\tilde{F}_{n_i}\circ
  \cdots \circ \tilde{F}_{n_1 +1}(z), \Phi_{n_i} \circ
\tilde{F}_{n_i}\circ \cdots  \circ
\tilde{F}_{n_1+1}(z))\leq \frac{1}{2^{i+1}}
\end{equation}
for all $z \in S_i$.
\begin{equation}\label{approxfixx}
\begin{split}
||D(\Phi_{n_{i+1}}& \circ \tilde{F}_{n_{i+1}}\cdots
\circ\tilde{F}_{n_{i}+1}\circ\tilde{F}_{n_i}\circ \cdots
 \circ
\tilde{F}_{n_1+1})(z)(v)||_{g_1}\\
&-||D(\Phi_{n_i} \circ \tilde{F}_{n_i}\circ \cdots  \circ
\tilde{F}_{n_1 +1})(z)(v)||_{g_1} \leq \frac{ b_i
}{2^{i+1}}
\end{split}
\end{equation}
for all $z \in S_i$ and Euclidean unit vectors $v$, where
the sequence $b_i$ is positive, decreases, and satisfies
\begin{equation}\label{approxfixxx}
||D\Phi_{n_i} \circ \tilde{F}_{n_i}\circ \cdot \cdot \circ
\tilde{F}_{2}(z)(v)||_{g_1}\geq b_i
\end{equation}
for all $z \in S_i$ and Euclidean unit vectors $v$. Here
$$
S_i=\lf(\tilde F_{n_i}\circ\cdots\circ \tilde
F_2\ri)^{-1}(D(\rho_3)).
$$
\end{lem}

\begin{cor}\label{fixc1}
Let $S_i$ be as above.  Then $S_i$ is an increasing
sequence of open sets in $\C^n$.
\end{cor}
\begin{proof}
from (\ref{fixe8}) we have
\begin{equation}\label{fixe10}
\tilde{F}_{n_{i+1}}\circ \cdots  \circ \tilde{F}_{n_i
+1}\circ \tilde{F}_{n_{i}}\circ \cdots  \circ
\tilde{F}_{2}(S_i)
 =\tilde{F}_{n_{i+1}}\circ \cdots \circ \tilde{F}_{n_i +1}(D(\rho_3))
 \subset D(\rho_3)
 \end{equation}
and thus
\begin{equation}\label{fixe11}
S_i\subset \tilde{F}_2^{-1}\circ \cdots \circ
\tilde{F}_{n_{i+1}}^{-1}(D(\rho_3)=S_{i+1}
 \end{equation}
\end{proof}
\begin{defn}\label{s3d1}
Let the sequences $S_i$ and $n_i$ be as above.  Let
$$
\Omega=\bigcup_{i=2}^\infty S_i
$$
\end{defn}

\begin{cor}\label{s3c1}
$\Omega$ is pseudoconvex and is homeomorphic to
$\mathbb{R}^{2n}$.
\end{cor}
\begin{proof}
Since each $S_i$ is pseduoconvex in $\C^n$, $\Omega$ is
pseudoconvex. Since each $S_i$ is homeomorphic to the unit
ball in $\mathbb{R}^{2n}$, $\Omega$ is also homeomorphic to
$\mathbb{R}^{2n}$ by \cite{B}.  This completes the proof of
the corollary.
\end{proof}
We now begin to use the maps $\tilde{F}_{i}$ to construct a map
from $\Omega$ onto $M$. We need the following lemma.
\begin{lem}\label{exhaustion1}
Let $M, g_i, p$ be as in Theorem \ref{main}.  Then for all $\e>0$,
$\bigcup_i B_i(\e)=M$, where $B_i(\e)=B_i(p,\e)$.
\end{lem}
\begin{proof}
Let $0<3\e<r_0$, where $r_0$ is as in \textbf{(a3)}.
 Obviously, $B_1(\e)\subset
\bigcup_{i}B_i(\e)$. We claim that if $B_1(k\e)\subset
\bigcup_{i}B_i(\e)$, $k\ge 1$, then
 $B_1((k+1)\e)\subset \bigcup_{i}B_i(\e)$.

 Suppose $B_1(k\e)\subset
\bigcup_{i}B_i(\e)$, then $\ol {B_1(k\e-\frac12\e)}\subset
B_{i}(\e)$ provided $i$ is large enough. Hence
$B_1((k+1)\e)\subset
 B_{i}(\e+\frac32\e)\subset B_{i}(3\e)$ for $i$ large
 enough by \textbf{(a1)}. By \textbf{(a3)}, we can find $i$
 such that  $B_1((k+1)\e)\subset B_i(\e)$. This completes
 the proof of the lemma.
\end{proof}
\begin{lem}\label{fixl2}
Let $\Gamma_i:=\Phi_{n_i}\circ \tilde{F}_{n_{i}}\circ \cdot
\cdot \circ \tilde{F}_{2}$.  Then the following map $\Psi :
\Omega \to M$ is well defined.
\begin{equation}\label{s3e4}
\Psi(z)=\lim_{i\to \infty}\Gamma_i (z).
\end{equation}
\end{lem}
\begin{proof}
This follows from (\ref{approxfix}) in Lemma
\ref{fixc2}, Corollary \ref{fixc1} and the definition of
the maps $\Gamma_i$.
\end{proof}

\begin{lem}\label{s3l2}
$\Psi$ is a local biholomorphism and onto.
\end{lem}
\begin{proof}

 By
 property (iv) of the maps $\Phi_i$, and the fact
 that $\bigcup_iB_i(p,\e)=M$ for all $\e$, given any $R>0$ we
 can find $n_i$ such that
\begin{equation}\label{s3e11}
B_1 (p, R) \subset B_{n_i} (p, \e) \subset
\Phi_{n_i}(D(C_1\e)),
\end{equation}
for some constant $C_1^2$ is the constant in Lemma
\ref{s1-L1}(iv). Here we have used  Corollary
\ref{cor-to-s1-L1} provided  $C_1\e<\rho_3$.
Choose such an $\e$.  Then
 $$\Gamma_i(S_i)=\Phi_{n_i}\circ\tilde{F}_{n_i}\circ \cdots \circ \tilde{F}_{2}(S_i)
 =\Phi_{n_i}(D(\rho_3)\supset B_1(p,R).
 $$
 Thus by (\ref{approxfix}) and the fact that the
$S_i$'s are increasing it follows that
\begin{equation}\label{s3e12}
B_1 (p,R-1) \subset  \Gamma_j (S_i)
\end{equation}
for all $j\geq i$.  From the definition of the map $\Psi$,
we see that
\begin{equation}\label{s3e13}
B_1 (p, R-1) \subset \Psi(\Omega).
\end{equation}
Hence $\Psi(\Omega)=M$.

 We now show that $\Psi$ is a local biholomorphism.
  Observe that $\Omega$ is open and $\Psi$ is a holomorphic map.
Now to show $\Psi$ is a local biholomorphism on $\Omega$,
it will be sufficient to show it is a local biholomorphism
on the sets $S_i$ for each $i$.   Fix some $i$.  Then by
(\ref{approxfixx}) and the fact that the $b_i$'s are
decreasing,
  \begin{equation}\label{approx4}
 ||D(\Gamma_j)(z)(v)||_{g_1} \geq  b_i - \frac{b_i}{2}
\end{equation}
  for all $j\geq i$, $z \in S_i$ and all unit
  vectors $v$ at $z$ .  Thus by the definition of $\Psi$,  (\ref{approx4}) implies $\Psi$ is a local biholomorphism on
  $S_i$.   Noting that $i$ is arbitrary, this completes the proof of the Lemma.

\end{proof}

\section{Proof of Theorem \ref{main}}

Let $M$ and $g_i$ satisfy \textbf{(a1)--(a3)}.  And let $\Psi$ be the map constructed in the previous section.
 If we take $\pi:\wh M\to M$ to be a universal holomorphic
covering of $M$ and let $\wh g_i=\pi^*(g_i)$, then $(\wh
M,\wh g_i)$ will still satisfy   \textbf{(a1)--(a3)}. Thus to prove Theorem \ref{main}, it will be sufficient to prove that $\Psi$  is injective assuming that $M$ is simply connected.  Before we
prove this, let us first prove the following:
\begin{lem}\label{homotopy} Let $\alpha(s)$, $0\le s\le 1$ be a smooth curve in
$M$. Then there exists $\e>0$ such that if $\beta(s)$ is
another smooth curve $M$ with same end points as
$\alpha(s)$ such that $d_1(\alpha(s),\gamma(s))<\e$ for all
$s$,  then there is a smooth homotopy $\gamma(s,\tau)$ with
end points fixed such that $\gamma(s,0)=\alpha(s)$ and
$\gamma(s,1)=\beta(s)$.  Moreover, there is a constant $L$
depending only on $(M,g_1)$, $\max_{0\le 1\le
1}\{|\alpha'(s)|_{g_1}+|\beta'(s)|_{g_1}\}$, such that the
length of $\gamma(\cdot,\tau)$ with respect to $g_1$ is
bounded above by $L$.
\end{lem}
\begin{proof} In the following, all lengths on $M$ will be computed with respect to the metric
 $g_0$. Let $\alpha(s)$ be given. Then there is $R>0$ such
that $\alpha\subset B_1(p,R/2)$. First let $\e>0$ be the
lower bound for the injectivity radius of $B_1(p, R)$.
Suppose $\beta(s)$ is another smooth curve on $M$ with same
end points as $\alpha(s)$. Then there is a smooth homotopy
$\gamma(s,\tau)$ such that $\gamma(s,\tau)$, for $0\le\tau
\le 1$  is the minimal geodesic from $\alpha(s)$ to
$\beta(s)$. Then for each $s$, $J=\gamma_s$ is a Jacobi
field along the geodesic $\gamma(s,\tau)$ for $0\le\tau \le
1$, with boundary value $J(0)= \alpha'(s)$ and
$J'(1)=\beta'(s)$. With respect to an orthonormal frame
$\{e_i\}$ parallel along $\gamma(s,\tau)$, $0\le\tau\le 1$,
the components $y_i$ of $J$ satisfies
$$
\lf[\begin{array}{crcl}
 y_1''\\ \vdots\\
y_{2n}''
\end{array}\ri]=A\lf[
\begin{array}{crcl}
y_1 \\
\vdots\\ y_{2n}\end{array} \ri]
$$
where $A_{ij}=\langle
R(\gamma_\tau,e_i)\gamma_\tau,e_j\rangle$. Here   $'$ means
derivatives with respect to $\tau$. Note that
$|\gamma_\tau|\le \e$ and the curvature is bounded from
below, we have
$$
\lf(\sum_{i}y_i^2\ri)''\ge -C_1\e^2\sum_{i}y_i^2
$$
for some constant $C_1>0$ depending only on the lower bound
of the curvature and $n$. Hence if $\e>0$ is small enough
depending only on the curvature, we can compare
$\sum_{i}y_i^2 $ with the solution $f$ of $f''=-C_1\e^2f$
with the same boundary value as $\sum_{i}y_i^2 $. Hence
$$
|\gamma_s|^2=|J|^2= \sum_{i}y_i^2  \le C_2
$$
for some $C_2$ depending only on $g_0$ and $\max_{0\le 1\le
1}\{|\alpha'(s)| +|\beta'(s)| \}$.
\end{proof}

 We now complete the proof of Theorem \ref{main} by proving
the following:
\begin{lem}\label{s33l4}
If $M$ is simply connected, then $\Psi$ is injective.
\end{lem}

\begin{proof}
Suppose the lemma is false.  Then there are distinct points
$z_1, z_2 \in \Omega$ such that $\Psi(z_1)=\Psi(z_2)=q$.
Let $\tilde{\gamma}(s)$ be a smooth curve in $\Omega$ for
$s\in [0, 1]$, joining $z_1$ to $z_2$ parametrized
proportional to arc length with respect to the Euclidean
metric, and let $\gamma(s)=\Psi\circ \tilde{\gamma}(s)$.
Then $\gamma(0)=\gamma(1)=q$.  Let $\gamma(s, \tau)$ be a
smooth homotopy of $\gamma$ for $(s, \tau)\in [0,
1]\times[0, 1]$ such that $\gamma(s, 0)=\gamma(s)$,
$\gamma(s, 1)=q$ for all $s\in [0, 1]$, and
$\gamma(0,\tau)=\gamma(1,\tau)=q$ for all $\tau\in [0,1]$.
Let $L_1=\max\{l(\gamma(\cdot,t)|\ \tau\in [0,1]\}$, where
$l(\gamma(\cdot,\tau))$ is the length of
$\gamma(\cdot,\tau)$ with respect to $g_1$.

Let $R>0$ be fixed, such that $\gamma(s,\tau)\in B_1(p,R)$
for all $0\le s,\tau\le1$.

By (\ref{s3e12}) and the fact that $S_i\subset S_{i+1}$ for
all $i$, there exists $i_0$ such that
\begin{equation}\label{inj1}
    B_1(p,R)\subset \Gamma_j(S_i)
\end{equation}
for all $j\ge i\ge i_0$, and that $\wt\gamma\subset
S_{i_0}$.

Since $\Psi$ is a local biholomorphism, it is easy to see
that for any $a>0$ there is $ b>0$ such that for all $i$
large enough, $\Gamma_i(D(z_k,a))\supset
B_1(\Gamma_i(z_k),b)$, $k=1, 2$. Since $\Gamma_i(z_k)\to
\Psi(z_k)=q$,  by choosing an even larger $ i_0$,  for all
$i\ge i_0$ there exist $\zeta_{1,i}\neq \zeta_{2,i}\in
S_{i_0}$ such that
$\Gamma_i(\zeta_{1,i})=\Gamma_i(\zeta_{2,i})=q$  and that
$\zeta_{1,i}\to z_1$ and $\zeta_{2,i}\to z_2$. Now for $i$
large enough, we can join $\zeta_{1,i}$ to $\zeta_{2,i}$ by
first joining $\zeta_{1,i}$ to $z_1$, then $z_1$ to $z_2$
along $\wt\gamma$, then $z_2$ to $\zeta_{2,i}$. Let us
denote this curve by $\wt\gamma_i(s)$, $s\in [0, 1]$
parametrized proportional to arc length. We may assume
$\wt\gamma_i(s)$ is smooth,
 $\wt\gamma_i(s)\subset K\subset S_{i_0}$ for some compact set $K$, and $|\wt\gamma_i'|\le C_1$ for some constant
independent of $i$ for all $i\ge i_0$.   Moreover, we have
$|\wt\gamma(s)-\wt\gamma_i(s)|\to0$ uniformly over $s$ as
$i\to\infty$. Since $\Psi$ is a local biholomorphism, there
is a constant $C_2$ independent of $i$ such that if
$\gamma_i=\Psi\circ \wt\gamma_i$, then
\begin{equation}\label{s5e1}
| \gamma_i'(s)|_{g_1}\le C_2.
\end{equation}
For the curve $\gamma(s)$, let $\e$ be as in Lemma
\ref{homotopy}. Since $\Gamma_i$ converge to $\Psi$
uniformly on compact sets together with first derivatives,
if $i_0$ is chosen large enough, then  the following are
true:
\begin{enumerate}
    \item [(i)] $d_1(\gamma(s),\gamma_i(s))<\frac\e2$
    \item [(ii)]
    $d_1(\Gamma_i\circ\wt\gamma_j(s),\gamma_j(s))
    =d_1(\Gamma_i\circ\wt\gamma_j(s),\Psi\circ\wt\gamma_j(s))<\frac\e2$
    \item [(iii)] $|(\Gamma_i\circ\wt\gamma_j)'(s)|_{g_1}\le
    |(\Psi\circ\wt\gamma_j)'(s)|_{g_1}+C_2=|
    \gamma_i'(s)|_{g_1}+ C_2\le 2C_2$
\end{enumerate}
for $i, j\ge i_0$.

By (i) and (ii), we have:
$$
d_1(\gamma(s),\Gamma_i\circ\wt\gamma_i(s))<\e
$$
for all $i\ge i_0$.  Thus by Lemma \ref{homotopy} and
(\ref{s5e1}), for each $i\ge i_0$ we can find a homotopy
which deforms $\gamma(s)$ to
$\Gamma_i\circ\tilde{\gamma}_i(s)$, with end points fixed,
so that each curve in the homotopy has length (with respect
to $g_1$) bounded by some constant $L$ independent of $i$.

 Now  let  $\rho_1$ be  as in Lemma \ref{s1-L2} corresponding to $\rho=r$.  Then
 we can choose $ i\ge i_0$ large enough but fixed, such that
 $B_1(p, L+L_1+R+1)\subset \Phi_{n_i}(D(\rho_3))$, and any curve $\beta$
  in the above homotopies is in $B_1(p,L+L_1+R+1)$ and satisfies $L_i(\beta)\le
 1/(L+L_1+R+1)\rho_3$. Here we have used \textbf{(a3)}.

 Let
 $w_k=\tilde F_{n_i}\circ\cdots\circ \tilde F_2(\zeta_{k,i})$, $k=1, 2$.
 Then $w_1\neq w_2$. Note that
 $$
\tilde F_{n_i}\circ\cdots\circ \tilde F_2(S_{i_0})\subset
\tilde F_{n_i}\circ\cdots\circ \tilde F_2(S_i)\subset
D(\rho_3).
$$
and $\rho_3<\frac18r$, see Remark \ref{s5r1}.

 By Corollary \ref{cor-to-L3}, since the lift of
 $\Gamma_i\circ\wt\gamma_i(s)$ in the Lemma \ref{s1-L2} from $w_1$
 by $\Phi_{n_i}$ is
 $\tilde F_{n_i}\circ\cdots\circ \tilde F_2\circ\wt\gamma_i(s)$, the lift
 $\wt\sigma$
  of $\gamma(\cdot,1)$ satisfies
  $\wt\sigma(1)  =\tilde F_{n_i}\circ\cdots\circ \tilde
 F_2\circ\wt\gamma_i(1)=w_2$.
This is impossible because
$\Phi_{n_i}\circ\wt\sigma(s)=\gamma(s,1)$ is a constant
map, $\wt\sigma(0)=w_1\neq w_2=\wt\sigma(1)$ and
$\Phi_{n_i}$ is a local biholomorphism.

\end{proof}

\section{Proof of Theorem \ref{s1t1}}
In this section we prove Theorem \ref{s1t1} .
 We begin proving a general theorem on complete solutions to the \KR flow

\begin{equation}\label{krf}
\frac{\p g_{i\jbar}}{\p t}=-R_{i\jbar}
\end{equation}

\begin{thm}\label{s6t1}
Let $g(t)$ be a complete solution to (\ref{krf}) with
non-negative holomorphic bisectional curvature such that
$g(0)$ has bounded curvature. Fix some $p \in M$ and let
$\lambda_i(t)$ be the eigenvalues of $Rc(p, t)$ arranged in
increasing order. Then
$$t\lambda_k(t)$$ is nondecreasing in  $t$ for all $1\leq k \leq n$.
\end{thm}

\begin{proof} To prove the theorem we may assume again that $M$ is simply
connected and by the result of \cite{cao1}, we may further
assume that the Ricci curvature is positive for all $x\in
M$ and for all $t>0$.

Now let $k\geq 1$, and let $h(t)$ be any positive function
with $h'(t)>0$ for all $t$. We claim that for any $t_0$
there is $\e>0$ such that
$th(t)\lambda_k(t)<t_0h(t_0)\lambda_k(t_0)$ for all $t\in
(t_0-\e,t_0)$.  By taking $h(t)=1+\delta t$ with $\delta>0$
and then let $\delta\to0$, we see that the theorem will
follow from this claim which we now prove.

  For any $t$, let  $0<\lambda_1(t) \le \dots \le \lambda_n(t)$ be the eigenvalues of
$R_\ijb (p, t)$.  For any $\sigma>0$ let $E_\sigma(t)$ be the direct sum
of the corresponding eigenspaces with eigenvalues $\lambda<\sigma$.
 Now let $t_0$ be fixed and let $m\ge k$ be the largest integer such that $\lambda_j(t_0)=\lambda_k(t_0)$ for
$m\ge j\ge k$. Let $\sigma>0$ be such that
$\lambda_k(t_0)<\sigma<\lambda_{m+1}(t_0)$ if $m<n$ and
$\sigma>\lambda_k(t_0)$ if $m=n$. Then there exists $\e>0$
such that for all $t\in (t_0-\e,t_0+\e)$,
$\lambda_m(t)<\sigma<\lambda_{m+1}(t)$ if $m<n$ and
$\sigma>\lambda_n(t)$ if $m=n$. In any case, the orthogonal
projection $P_\sigma(t)$ onto $E_\sigma(t)$ is smooth in
$(t_0-\e,t_0+\e)$, see \cite[p.501]{CT3} for example.   For
any $t_1\in (t_0-\e,t_0+\e)$, let $v_1$ be an eigenvector
of $R_\ijb(t)$ corresponding to $\lambda_m(t_1)$ with
length 1. Let
$$v(t)=\frac{P_\sigma(t)v_1}{|P_\sigma(t)v_1|_t}.
$$
Note that for $t$ close to $t_1$, $ P_\sigma(t)v_1\neq0$.
In local holomorphic coordinates $z^i$, let $a(t)=R_\ijb
v^i\bar v^j$ where $v(t)=v^i(t)\frac{\p}{\p z^i}$. Note
that $P_\sigma (t_1) (v_1)=v_1$ and thus
$a(t_1)=\lambda_m(t_1)$.  Also note that we have $a(t)\leq \lambda_m (t)$ for all $t\in (t_0-\e,t_0+\e)$ where $a(t)$ is defined.  Now note that
\begin{equation}
0=\frac{d}{dt} \langle v(t),v(t)\rangle_t   =
-R_{i\jbar}(p, t)v^i v^{\bar{j}} + 2Re (g_{i \jbar}
  \frac{d v^i}{dt} v^{\bar{j}} ).
\end{equation}
Hence at $t_1$ we have
\begin{equation} \label{s6e1}
   Re (g_{i \jbar}
  \frac{d v^i}{dt} v^{\bar{j}} )=\frac{a(t_1)}2=\frac{\lambda_m(t_1)}2
\end{equation}
 By the Harnack Inequality in \cite{Cao97} we have
\begin{equation}\label{harnack}
\frac{\partial R_{i\jbar}}{\partial t} + g^{k \bar{l}} R_{i
\bar{l}} R_{k \bar{j}} + \frac{R_{i\jbar}}{t} \geq 0.
\end{equation}
for all $t$. Thus at $t_1$ we have

  \begin{equation}\label{mde1}
 \begin{split}
0 &\leq \frac{\partial R_{i\jbar}}{\partial t}v^i v^{\jbar} +
g^{k \bar{l}} R_{i \bar{l}} R_{k \bar{j}} v^i v^{\jbar}+
\frac{R_{i\jbar}}{t_1}v^i v^{\jbar}\\
&=\frac{d}{dt}( R_{i\jbar}v^i v^{\jbar} )- 2Re
\lf(R_{i\jbar}(\frac{d}{dt}v^i) v^{\jbar}\ri )
 + g^{k \bar{l}} R_{i \bar{l}} R_{k \bar{j}} v^i v^{\jbar}
  + \frac{R_{i\jbar}}{t}v^i v^{\jbar}\\
&=\frac{d}{dt}( R_{i\jbar}v^i v^{\jbar} )-\lambda_m^2(t_1)
+ \lambda_m^2(t_1)
+\frac{R_{i\jbar}}{t}v^i v^{\jbar}\\
&=\frac{d}{dt}( R_{i\jbar}v^i v^{\jbar} ) +\frac{R_{i\jbar}}{t_1}v^i v^{\jbar}\\
&=\frac{d}{dt}a+ \frac{a}{t_1}.
  \end{split}
  \end{equation}
  where the third  equality follows from writing the expressions
  in a holomorphic coordinate $z^i$
  such that $\frac{\p}{\p z^i}$ form  a basis of eigenvectors
  of $Rc(p, T)$ with $v_1=\frac{\p}{\p z^1}$ at $p$,
   and (\ref{s6e1}).
Since $a(t_1)>0$, $h(t)>0$ and $h'(t)>0$, we conclude that
$$\frac{d}{dt}\lf(th(t)a(t)\ri)>0
$$
at $t_1$ and hence   $th(t)a(t)$ is increasing in $t$ for
$t\in (t_1-\e_1,t_1+\e_1)$ for some $\e_1>0$. Hence
\begin{equation} \begin{split}
   t_1h(t_1)\lambda_m(t_1)&=t_1h(t_1)a(t_1)\\
        &<th(t)a(t)\\
        &\le th(t)\lambda_m(t)
\end{split}
\end{equation}
for all $t_1<t<t_1+\e_1$.  As $t_1$ was chosen arbitrarily
in $(t_0-\e,t_0+\e)$,  we conclude that $th(t)\lambda_m(t)$
is increasing in $(t_0-\e,t_0+\e)$. In particular,
\begin{equation} \begin{split}
   t_0h(t_0)\lambda_k(t_0)&=t_0h(t_0)\lambda_m(t_0)\\
        &> th(t)\lambda_m(t)\\
        &\ge th(t)\lambda_k(t)
\end{split}
\end{equation}
for all $t\in (t_0-\e,t_0)$. This proves the claim and the
theorem.
\end{proof}
\begin{proof}[Proof of Theorem \ref{s1t1}]
We begin by observing that if $M$ has positive holomorphic
bisectional curvature and is simply connected near
infinity, then it is actually simply connected. Indeed, if
$M$ were not simply connected there would exist a
nontrivial minimizer of a free homotopy class.  This
however is impossible by the fact that the bisectional
curvature is positive, and an argument as in the proof of
Sygne theorem.  By   the remarks at the beginning of $\S$
5, we may assume that $M$ is simply connected in Theorem
\ref{s1t1} and by \cite{cao1} we may also assume that the
Ricci curvature is positive in spacetime.

Now by the long time existence results in \cite{Sh2} we
know that under the hypothesis of Theorem 1.1, (\ref{krf})
has a long time solution $g(t)$ with uniformly bounded
non-negative holomorphic bisectional curvature together
with the covariant derivatives of the curvature tensor.
Since  $Rc>0$, Theorem \ref{s6t1} implies that given any
compact set $\Omega$  we can find $C>0$ such that $Rc(t)
\geq \frac{C}{t} g(t)$ on $\Omega$ for all $t$. From this,
(\ref{krf}) and recalling that the curvature of $g(t)$ is
uniformly bounded on $[0, \infty)\times M$, it is not hard
to see show that the sequence of metrics $g_i=g(i)$ on $M$
satisfies the hypothesis of Theorem 1.2 and thus Theorem
1.1 follows.
\end{proof}

\bibliographystyle{amsplain}

\begin{thebibliography}{10}

\bibitem{A} Anderson, E. and Lempert, L., {\sl On the group of holomorphic automorphisms of $\C^n$},  Invent. Math.  110  (1992),  no. 2, 371--388.
\bibitem{B}Brown, M., {\sl The monotone union of open $n$-cells is an open
$n$-cell},  Proc. Amer. Math. Soc. \textbf{12} (1961), 812--814.

\bibitem{Cao92}  Cao, H.-D.,
{\sl On Harnack's inequality for the K\"ahler-Ricci flow},
Invent. Math. \textbf{109} (1992),
  247--263.


\bibitem{Cao97} Cao, H.-D., {\sl Limits of solutions to the
  K\"ahler-Ricci flow}, J. Differential Geom.\ \textbf{45} (1997),
  257--272.

   \bibitem{cao1} Cao, H.-D., {\sl On Dimension reduction in the K\"ahler-Ricci flow},  Comm. Anal. Geom. \textbf{12} (2004),
  305--320.


\bibitem{CT} Chau, A. and Tam, L.-F., {\sl Gradient
\KRS and a uniformization conjecture}, arXiv eprint 2002.
arXiv:math.DG/0310198.

\bibitem{CT2} Chau, A. and Tam, L.-F., {\sl A note on the uniformization
of gradient K\"ahler-Ricci solitons }, Math. Res. Lett.\textbf{
12} (2005),   19--21..

 \bibitem{CT3} Chau, A. and Tam, L.-F.,
 {\sl On the complex structure of K\"{a}hler manifolds with non-negative
 curvature}, J. Differential Geom. \textbf{73} (2006),   491--530.

\bibitem{CT4} Chau, A. and Tam, L.-F.,{\sl Non-negatively
curved K\"ahler manifolds with average quadratic curvature decay},
to appear in Communications in Analysis and Geometry.

\bibitem{Chen} Chen, B.L. and Zhu, X.P.,
 {\sl On complete noncompact K\"ahler manifolds with positive
 bisectional curvature},   Math. Ann. \textbf{327} (2003), 1--23.

\bibitem{CTZ} Chen, B.L., Tang, S.H. and Zhu, X.P.,
 {\sl A Uniformization Theorem Of Complete Noncompact
       K\"{a}hler Surfaces With Positive Bisectional Curvature},  J. Differential Geom.
       \textbf{ 67} (2004),  519--570.

\bibitem{CZ2} Chen, B.L. and Zhu, X.P.,
 {\sl Positively Curved Complete Noncompact K\"{a}hler
       Manifolds},  arXiv eprint 2002. arXiv:math.DG/0211373.

\bibitem{CZ3} Chen, B.L. and Zhu, X.P.,
 {\sl Volume Growth and Curvature Decay of Positively
       Curved K\"{a}hler manifolds},   Q. J. Pure Appl. Math.  \textbf{1 }(2005),  68--108.

\bibitem{GW1} Greene, R. E. and Wu, H.,{\sl Analysis on noncompact K\"ahler manifolds}, Proc. Sympos. Pure Math., \textbf{30} Part 2 (1977), 69-100.

\bibitem{GW76} Greene, R. E. and Wu, H., {\sl $C^\infty$
convex function and the manifolds of positive curvature},
Acta. Math., \textbf{137} (1976),   209--245.

\bibitem{GM} Gromoll, D. and Meyer, W., {\sl On complete open manifolds of positive curvature},  Ann. of Math. \textbf{90} (1969), 75--90.

\bibitem{Hamiton95} Hamilton, R. S.,
 {\sl Formation of Singularities in the Ricci Flow},
 Surveys in differential geometry, Vol. II (1995),  7--136.



\bibitem{M1}  Mok, N., {\sl An embedding theorem of complete K\"a
  manifolds of positive bisectional curvature onto affine algebraic varieties},
Bull. Soc. Math. France. \textbf{112} (1984), 179--258.

\bibitem{M2}   Mok, N., {\sl An embedding theorem of complex K\"ahler
  manifolds of positive Ricci curvature onto quasi-projective varieties },
Math. Ann. \textbf{286} (1990), no.1-3,  373--408.

\bibitem{MSY}   Mok, N.,  Siu, Y.-T. and Yau, S.-T., {\sl The Poincar\'e-Lelong
  equation on complete K\"ahler manifolds },
Comp. Math.,\textbf{44},(1981), 183--218.

\bibitem{Ni} Ni, L., {\sl Vanishing theorems on complete K\"ahler manifolds and their applications},
J. Differential Geom. \textbf{50} (1998), no  89--122.

\bibitem{Ni2} Ni, L., {\sl Ancient solutions to K\"ahler-Ricci flow}, Math. Res. Lett.
\textbf{12} (2005),  633--653.

\bibitem{Ni0} Ni, L., {\sl A new Li-Yau-Hamilton estimate for
  K\"ahler-Ricci flow }, arXiv eprint 2005. arXiv:math.DG/0502495.


\bibitem{NST} Ni, L., Shi, Y.-G. and  Tam, L.-F.
{\sl Poisson equation, Poincar\'e-Lelong  equation and curvature
decay on complete K\"ahler manifolds},  J. Differential Geom.
\textbf{57} (2001),  339--388.


\bibitem{NT} Ni, L. and Tam, L.-F.,
 {\sl \KRF and the Poincar\'e-Lelong equation},  Comm. Anal. Geom.
 \textbf{12} (2004), 111--141.

\bibitem{NT2} Ni, L. and Tam, L.-F.,
 {\sl Plurisubharmonic functions and the structure of complete K\"ahler
 manifolds with nonnegative curvature},  J. Differential Geom. \textbf{64} (2003),   457--524.


\bibitem{Sh0} Shi, W.-X., {\sl Ricci deformation of the metric on complete noncompact Riemannian manifolds}, J. of Differential Geometry  \textbf{30} (1989), 223-301.


\bibitem{Sh} Shi, W.-X.,
 {\sl Ricci deformation of the metric on complete noncompact \K manifolds},
  PhD thesis, Harvard University, 1990.

\bibitem{Sh1} Shi, W.-X.,{\sl Complete noncompact K\"ahler manifolds with positive holomorphic bisectional curvature}, Bull. Amer. Math. Soc. (N. S.) \textbf{23}
(1990), 437--400.

\bibitem{Sh2} Shi, W.-X.,
 {\sl Ricci Flow and the uniformization on
complete non compact \K manifolds},
  J. of Differential Geometry  \textbf{45} (1997),  94-220.

 \bibitem{Siu}    Siu, Y.-T., {\sl
Pseudoconvexity and the problem of Levi}, Bull. Amer. Math. Soc.
\textbf{84} (1978), 481--512.
\bibitem{Y1} Yau, S.-T., {\sl Some Function-Theoretic Properties of Complete Riemannian Manifold and Their Applications to
Geometry}, Indiana University Mathematics Journal.
 Vol. 25, No. 7 (1976), 659--670.

\bibitem{Y} Yau, S.-T., {\sl A review of complex differential geometry},
 Proc. Sympos. Pure Math., \textbf{52} Part 2 (1991), 619--625.
\bibitem{W1}   Wu, H., {\sl  An elementary methods in the
study of nonnegative curvature},  Acta. Math. \textbf{ 142 } (
1979),   57--78.

\bibitem{Z} Zhu, X.P.,
 {\sl The Ricci Flow on Complete Noncompact K\"{a}hler
       Manifolds},  arXiv eprint 2002. arXiv:math.DG/0211375.
\end{thebibliography}

\end{document}